\documentclass[12pt,reqno]{amsart}

\usepackage{amssymb}
\usepackage{fancyhdr}

\usepackage{eucal}

\input{xypic}

\def\cL{{\mathcal L}}

\def\cR{{\mathcal R}}

\def\BZ{{\mathbb Z}}

\def\sA{\mbox{\sf A}}
\def\sB{\mbox{\sf B}}

\def\sD{\mbox{\sf D}}

\def\sF{\mbox{\sf F}}
\def\sG{\mbox{\sf G}}

\def\sI{\mbox{\sf I}}

\def\sP{\mbox{\sf P}}
\def\sQ{\mbox{\sf Q}}

\def\b{\operatorname{b}}

\def\Ext{\operatorname{Ext}}

\def\fd{\operatorname{fd}}

\def\H{\operatorname{H}}

\def\Hom{\operatorname{Hom}}
\def\id{\operatorname{id}}

\def\inf{\operatorname{inf}}

\def\LTensor{{\otimes}^{\operatorname{L}}}

\def\RHom{\operatorname{RHom}}

\def\sup{\operatorname{sup}}

\def\ten{\otimes}

\def\Tor{\operatorname{Tor}}

\numberwithin{equation}{part}

\newtheorem{Lemma}{Lemma}[section]
\newtheorem{Theorem}[Lemma]{Theorem}
\newtheorem{Proposition}[Lemma]{Proposition}
\newtheorem{Corollary}[Lemma]{Corollary}

\theoremstyle{definition}
\newtheorem{Definition}[Lemma]{Definition}

\newtheorem{Remark}[Lemma]{Remark}

\newtheorem{Observation}[Lemma]{Observation}
\newtheorem{Example}[Lemma]{Example}

\def\skewtimes{\ltimes}

\def\Gpd{\operatorname{Gpd}}
\def\Gfd{\operatorname{Gfd}}
\def\Gid{\operatorname{Gid}}

\def\CGpd{\mbox{$\operatorname{Gpd}_{A \ltimes C}$}}
\def\CGfd{\mbox{$\operatorname{Gfd}_{A \ltimes C}$}}
\def\CGid{\mbox{$\operatorname{Gid}_{A \ltimes C}$}}

\def\AGpd{\mbox{$\operatorname{Gpd}_{A \ltimes A}$}}
\def\AGfd{\mbox{$\operatorname{Gfd}_{A \ltimes A}$}}
\def\AGid{\mbox{$\operatorname{Gid}_{A \ltimes A}$}}

\def\GdimC{\mbox{$G$--$\operatorname{dim}_C$}}
\def\GdimAxC{\mbox{$G$--$\operatorname{dim}_{A \ltimes C}$}}

\def\GIdim{\mbox{$C\textnormal{-Gid}_A$}}
\def\GPdim{\mbox{$C\textnormal{-Gpd}_A$}}
\def\GFdim{\mbox{$C\textnormal{-Gfd}_A$}}

\def\AGIdim{\mbox{$A\textnormal{-Gid}_A$}}
\def\AGPdim{\mbox{$A\textnormal{-Gpd}_A$}}
\def\AGFdim{\mbox{$A\textnormal{-Gfd}_A$}}

\def\ProperGIdim{\mbox{$C\textnormal{-}\mathsf{Gid}_A$}}
\def\ProperGPdim{\mbox{$C\textnormal{-}\mathsf{Gpd}_A$}}
\def\ProperGFdim{\mbox{$C\textnormal{-}\mathsf{Gfd}_A$}}

\begin{document}

\title[Semi-dualizing modules and Gorenstein dimensions]
{Semi-dualizing modules and related Gorenstein homological dimensions}

\author{Henrik Holm \ \ }

\address{ {\flushleft Matematisk} Afdeling, Universitetsparken 5,
  DK--2100 K\o benhavn \O, Denmark, e-mail address:
  \textnormal{\texttt{holm@math.ku.dk}} }



\author{\ \ Peter J\o rgensen} 

\address{ {\flushleft Department} of Pure Mathematics, University of
  Leeds, Leeds LS2 9JT, United Kingdom, e-mail address:
  \textnormal{\texttt{popjoerg@maths.leeds.ac.uk}}, homepage:
  \textnormal{\texttt{http://www.maths.leeds.ac.uk/\raisebox{-0.6ex}{\~}{}popjoerg}} }   




\begin{abstract} 
  A semi-dualizing module over a commutative noetherian ring $A$ is a
  finitely generated module $C$ with $\RHom_A(C,C) \simeq A$ in the
  derived category $\sD(A)$. 
  
  \medskip
  \noindent
  We show how each such module gives rise to three new homological
  dimensions which we call $C$--Gorenstein projective, $C$--Gorenstein
  injective, and $C$--Gorenstein flat dimension, and investigate the
  properties of these dimensions.
\end{abstract}

\maketitle
\setcounter{section}{0}


\section*{Introduction}

\noindent
It is by now a well-established fact that over any associative ring
$A$, there exists a Gorenstein injective, Gorenstien projective and
Gorenstein flat dimension defined for complexes of $A$--modules. These
are usually denoted $\Gid_A(-)$, $\Gpd_A(-)$ and $\Gfd_A(-)$,
respectively. Some references are \cite{Wintherbook}, \cite{CFH},
\cite{HHGorensteinHomDim}, and \cite{Veliche}.

\medskip
\noindent
In this paper, we need to consider \emph{semi-dualizing} $A$--modules
$C$ (see Definition \ref{def:semi-dual}), and in order to make things
less technical, we only consider commutative and noetherian rings.

\medskip
\noindent
For any semi-dualizing module (in fact, complex) $C$ over $A$, and any
complex $Z$ with bounded and finitely generated homology, Christensen
\cite{Wintherpaper} introduced the dimension $\GdimC Z$, and developed
a satisfactory theory for this new invariant.

\medskip
\noindent
If $C$ is a semi-dualizing $A$--module and $M$ is any $A$--complex,
then we suggested in \cite{HHPJ2} the viewpoint that one should change
rings from $A$ to $A \ltimes C$ (the \emph{trivial extension} of $A$
by $C$; see Definition \ref{def:trivial_extension}), and then consider
the three ``ring changed'' Gorenstein dimensions:
\begin{displaymath}
  \CGid M \quad , \quad \CGpd M \quad \textnormal{ and } \quad \CGfd
  M.
\end{displaymath}
The usefulness of this viewpoint was demonstrated as it enabled us to
introduce three new \emph{Cohen-Macaulay dimensions}, which
characterize Cohen-Macaulay rings in a way one could hope for.

\medskip
\noindent
In this paper, we define for every semi-dualizing $A$--module $C$,
three new Gorenstein dimensions:
\begin{displaymath}
  \GIdim(-) \quad , \quad \GPdim(-) \quad \textnormal{ and } \quad
  \GFdim(-),
\end{displaymath}
which are called the \emph{$C$--Gorenstein injective},
\emph{$C$--Gorenstein projective} and \emph{$C$--Gorenstein flat
  dimension}, respectively (see Definition \ref{def:new_dimensions}).

\medskip
\noindent
It is worth pointing out that the, say, $C$--Gorenstein injective
dimension is defined in terms of resolutions consisting of so-called
\emph{$C$--Gorenstein injective $A$--modules} (see Definition
\ref{def:Cmodules}); and it does not involve a change of rings.  The
$C$--Gorenstein dimensions have at least five nice properties:
\smallskip
\begin{enumerate}
\item For complexes with bounded and finitely generated homology, our
  $\GPdim(-)$ agrees with Christensen's $\GdimC(-)$ (Proposition
  \ref{pro:GdimC}). \medskip
\item The three $C$--Gorenstein dimensions always agree with the
  ``ring changed'' dimensions $\CGid(-)$, $\CGpd(-)$ and $\CGfd(-)$,
  which were so important in \cite{HHPJ2} (Theorem
  \ref{thm:main_theorem}). \medskip
\item If $C=A$, the $C$--Gorenstein dimensions agree with the
  classical Gorenstein dimensions $\Gid_A(-)$, $\Gpd_A(-)$ and
  $\Gfd_A(-)$.
\end{enumerate}
If $A$ admits a dualizing complex $D$; cf.~\cite[Definition
(1.1)]{CFH}, then finiteness of the $C$--Gorenstein dimensions can be
interpreted in terms of \emph{Auslander and Bass categories} (see
Remark \ref{rmk:AB}):
\begin{enumerate}
\item[$(4)$] If we define $C^{\dagger} = \RHom_A(C,D)$, then for all
  (appropriately homologically bounded) $A$--complexes $M$ and $N$, we
  have the following implications (Theorem \ref{thm:AB}):
  \begin{align*}
    M \in \sA_{C^{\dagger}}(A) &\,\Leftrightarrow\, \GPdim M<\infty
    \,\Leftrightarrow\, \GFdim M<\infty; \\
    N \in \sB_{C^{\dagger}}(A) &\,\Leftrightarrow\, \GIdim N<\infty.
  \end{align*}
\end{enumerate}
This generalizes the main results in \cite[Theorems (4.3) and
(4.5)]{CFH}.

\medskip 
\noindent
Finally, each of the three $C$--Gorenstein dimensions has a related
\emph{proper} variant, giving us three additional dimensions
(Definitions \ref{dfn:ProperDimensions} and
\ref{dfn:SpecialProperDimensions}):
\begin{displaymath}
  \ProperGIdim(-) \quad , \quad \ProperGPdim(-) \quad \textnormal{ and
    } \quad \ProperGFdim(-).
\end{displaymath}
It turns out that the best one could hope for really happens, as we in
Theorems \ref{thm:proper_CGid}, \ref{thm:proper_CGpd} and
\ref{thm:proper_CGfd} prove:

\begin{enumerate}
\item[$(5)$] The proper $C$--Gorenstein dimensions (whenever these are
  defined) agree with the ordinary $C$--Gorenstein dimensions.
\end{enumerate}

\noindent
The paper is organized as follows: 

\medskip
\noindent
In Section \ref{sec:AxC} we have collected some fundamental facts
about the trivial extension $A \ltimes C$, which will be important
later on.  Section \ref{sec:CGdims} defines the three new
$C$--Gorenstein dimensions and proves how they are related to the
``ring changed'' Gorenstein dimensions over $A \ltimes C$.  Section
\ref{sec:Christensen} compares our $\GPdim(-)$ with Christensen's
$\GdimC(-)$.  In Section \ref{sec:AB} we interpret the $C$--Gorenstein
dimensions in terms of Auslander and Bass categories.  Finally,
Section \ref{sec:proper_dimensions} investigates the proper
$C$--Gorenstein dimensions.

\medskip
\noindent
{\bf Setup and notation.}\
  Throughout this paper, $A$ is a fixed commutative and noetherian
  ring with unit, and $C$ is a fixed semi-dualizing $A$--module;
  cf.~Definition \ref{def:semi-dual} below.

  \medskip
  \noindent
  We work within the derived category $\sD(A)$ of the category of
  $A$--modules; cf.~e.g.~\cite[Chapter I]{HartsResDual} and
  \cite[Chapter 10]{Weibel}; and complexes $M \in \sD(A)$ have
  differentials going to the right:
  \begin{displaymath}
    M = \xymatrix{\cdots \ar[r] &  M_{i+1} \ar[r]^-{\partial^M_{i+1}} &
      M_i \ar[r]^-{\partial^M_i} & M_{i-1} \ar[r] & \cdots}.
  \end{displaymath}
  We consistently use the hyper-homological notation from
  \cite[Appendix]{Wintherbook}, in particular we use $\RHom_A(-,-)$
  for the right derived Hom functor, and $-\LTensor_A-$ for the left
  derived tensor product functor.


\section{A few results about the trivial extension}
\label{sec:AxC}

\noindent
In this section we collect some fundamental results about the trivial
extension, which will be important later on.

\begin{Definition} \label{def:semi-dual}
  A finitely generated $A$--module $C$ with $\RHom_A(C,C) \simeq A$ in
  $\sD(A)$ is called \emph{semi-dualizing} ($C=A$ is such an example). 
\end{Definition}

\begin{Definition} \label{def:trivial_extension}
  If $C$ is any $A$--module, then the direct sum $A \oplus C$ can be
  equipped with the product:
  \begin{displaymath}
    (a,c)\cdot(a',c') \,=\, (aa',ac'+a'c).
  \end{displaymath}
  This turns $A \oplus C$ into a ring which is called the
  \emph{trivial extension} of $A$ by $C$ and denoted $A \skewtimes C$.

  \medskip
  \noindent
  There are canonical ring homomorphisms, $\!\!\xymatrix{A \ar@<0.5ex>[r]
    & A \ltimes C, \ar@<0.5ex>[l]}\!$ which enable us to view
  $A$--modules as $(A \ltimes C)$--modules, and vice versa. This will
  be done frequently.
\end{Definition}

\noindent
We import from \cite[Lemma 3.2]{HHPJ2} the following facts about
the interplay between the rings $A$ and $A \ltimes C$:

\begin{Lemma} \label{lem:extension_formulae}
  Let $A$ be a ring with a semi-dualizing module $C$.
  \begin{enumerate}
  \item There is an isomorphism in $\sD(A\ltimes C)$:
    \begin{displaymath}
      \RHom_A(A \skewtimes C,C) \cong A \skewtimes C.
    \end{displaymath}
  \item There is a natural equivalence of functors on $\sD(A)$:
    \begin{displaymath}
       \RHom_{A \skewtimes C}(-,A \skewtimes C) \simeq \RHom_A(-,C). 
    \end{displaymath}
  \item If $M$ is in $\sD(A)$ then the two biduality morphisms:
    \begin{displaymath}
    \begin{array}{c}
      M \longrightarrow \RHom_A(\RHom_A(M,C),C) 
      \ \text{ and } \medskip \\
      M \longrightarrow \RHom_{A \skewtimes C}(\RHom_{A \skewtimes
        C}(M,A \skewtimes C),A \skewtimes C) 
     \end{array}
     \end{displaymath}
    are equal.
  \item There is an isomorphism in $\sD(A \ltimes C)$:
    \begin{xxalignat}{3}
      &{\hphantom{\qed}}
      & \RHom_{A \skewtimes C}(A,A \skewtimes C) \,\cong\, C.
      & {}
      && \qed
    \end{xxalignat}
  \end{enumerate}
\end{Lemma}

\noindent
Furthermore, we have the next result \cite[Lemma 3.1]{HHPJ2} about
injective modules over $A$ and $A \ltimes C$:

\begin{Lemma} \label{lem:induced_injectives}
  The following two conclusions hold:
  \begin{enumerate}
  \item If $I$ is a (faithfully) injective $A$--module then $\Hom_A(A
    \skewtimes C,I)$ is a (faithfully) injective $(A \skewtimes
    C)$--module.
  \item Each injective $(A \skewtimes C)$--module is a direct summand
    in a module $\Hom_A(A \skewtimes C,I)$ where $I$ is some injective
    $A$--module. \hspace{\stretch{1}} $\qed$
  \end{enumerate}
\end{Lemma}

\noindent
Using the same methods, we obtain:

\begin{Lemma} \label{lem:induced_projectives}
  The following two conclusions hold:
  \begin{enumerate}
  \item If $P$ is a projective $A$--module then $(A \skewtimes
    C)\ten_AP$ is a projective $(A \skewtimes C)$--module.
  \item Each projective $(A \skewtimes C)$--module is a direct summand
    in a module $(A \skewtimes C)\ten_AP$ where $P$ is some projective
    $A$--module. \hspace{\stretch{1}} $\qed$
  \end{enumerate}
\end{Lemma}


\section{$C$--Gorenstein homological dimensions} \label{sec:CGdims}

\noindent
Let $M$ be an (appropriately homologically bounded) $A$--complex.  In
\cite{HHPJ2} we demonstrated the usefulness of changing rings from $A$
to $A \ltimes C$, and then considering the ``ring changed'' Gorenstein
dimensions:
\begin{displaymath}
  \CGid M \quad , \quad \CGpd M \quad \textnormal{ and } \quad \CGfd
  M.
\end{displaymath}
This point of view enabled us to introduce three \emph{Cohen-Macaulay
dimensions} which characterize Cohen-Macaulay local rings in a way one
could hope for. The next result is taken from \cite[Lemma 4.6]{HHPJ2}.


\begin{Proposition} \label{prop:dual}
  If $E$ is a faithfully injective $A$--module, and $M$ is any
  homologically right-bounded $A$--complex, then:
   \begin{xxalignat}{3}
      &{\hphantom{\qed}}
      &  \CGid  \Hom_A(M,E) \,=\, \CGfd  M.
      & {}
      && \qed
    \end{xxalignat}
\end{Proposition}


\begin{Lemma} \label{lem:calculation}
  Let $J$ be an injective $A$--module and $Q$ a projective
  $A$--module. Then we have a natural equivalence of functors on
  $\sD(A \ltimes C)$:
  \begin{enumerate}
  \item $\RHom_{A \ltimes C}(\Hom_A(A \ltimes C,J),-) \,\simeq\,
    \RHom_A(\Hom_A(C,J),-)$. \smallskip
  \item $\RHom_{A \ltimes C}(-,(A \ltimes C)\ten_AQ) \,\simeq\,
    \RHom_A(-,C\ten_AQ)$. 
  \end{enumerate}
\end{Lemma}

\begin{proof}
  (1) is \cite[Lemma 3.4]{HHPJ2}, and (2) is proved similarly.
\end{proof}


\begin{Corollary} \label{cor:Ext}
  For any $A$--module $M$, and integer $n$ we have:
  \begin{enumerate}
  \item $\Ext^n_A(\Hom_A(C,J),M)=0$ for all injective $A$--modules $J$
    if and only if $\Ext^n_{A \ltimes C}(U,M)=0$ for all injective $(A
    \ltimes C)$--modules $U$.
  \item $\Ext^n_A(M,C\ten_AP)=0$ for all projective $A$--modules $P$
    if and only if $\Ext^n_{A \ltimes C}(M,S)=0$ for all projective
    $(A \ltimes C)$--modules $S$.
  \end{enumerate}
\end{Corollary}

\begin{proof}
  (1) follows from Lemmas \ref{lem:calculation}(1) and
  \ref{lem:induced_injectives}, while (2) is a consequence of Lemmas 
  \ref{lem:calculation}(2) and \ref{lem:induced_projectives}.
\end{proof}

\noindent
We need to recall the next result from \cite[Lemma 4.1]{HHPJ2}. Its
proof uses, in fact, Lemmas \ref{lem:calculation}(1) and
\ref{lem:induced_injectives}.

\begin{Lemma} \label{lem:CGorInj/A<=GorInj/AxC}
  Let $M$ be an $A$--module which is Gorenstein injective over $A
  \skewtimes C$. Then there exists a short exact sequence of
  $A$--modules,
  \begin{displaymath}
    0 \rightarrow M^{\prime} \longrightarrow \Hom_A(C,I)
    \longrightarrow M \rightarrow 0, 
  \end{displaymath}
  where $I$ is injective over $A$ and $M^{\prime}$ is Gorenstein
  injective over $A \skewtimes C$.  Furthermore, the sequence stays
  exact if one applies to it the functor $\Hom_A(\Hom_A(C,J),-)$ for
  any injective $A$--module $J$. \hspace{\stretch{1}} $\qed$
\end{Lemma}

\noindent
``Dualizing'' the proof of Lemma \ref{lem:CGorInj/A<=GorInj/AxC}; this
time using Lemmas \ref{lem:calculation}(2) and
\ref{lem:induced_projectives}, we establish the next:

\begin{Lemma} \label{lem:CGorProj/A<=GorProj/AxC}
  Let $M$ be an $A$--module which is Gorenstein projective over $A
  \skewtimes C$. Then there exists a short exact sequence of
  $A$--modules,
  \begin{displaymath}
    0 \rightarrow M \longrightarrow C\ten_AP \longrightarrow M'
    \rightarrow 0,  
  \end{displaymath}
  where $P$ is projective over $A$ and $M^{\prime}$ is Gorenstein
  projective over $A \skewtimes C$.  Furthermore, the sequence stays
  exact if one applies to it the functor $\Hom_A(-,C\ten_A Q)$ for any
  projective $A$--module $Q$.  \hspace{\stretch{1}} $\qed$
\end{Lemma}

\noindent
The last result we will need to get started is \cite[Lemma
3.3]{HHPJ2}:

\begin{Lemma} \label{lem:Proj/A=>GorProj/AxC}
  The $A$--modules $A$ and $C$ are Gorenstein projective over $A
  \skewtimes C$. If $I$ is an injective $A$--module, then $\Hom_A(A,I)
  \cong I$ and $\Hom_A(C,I)$ are Gorenstein injective over $A \ltimes
  C$.  \hspace{\stretch{1}} $\qed$
\end{Lemma}

\noindent
Next, we introduce three new classes of modules:

\begin{Definition} \label{def:Cmodules}
  An $A$--module $M$ is called \emph{$C$--Gorenstein injective} if:
  \begin{enumerate}
  \item[(I1)] $\Ext_A^{\geqslant 1}(\Hom_A(C,I),M)=0$ for all injective
    $A$--modules $I$.
  \item[(I2)] There exist injective $A$--modules $I_0,I_1,\ldots$
    together with an exact sequence:
    \begin{displaymath}
       \cdots \to \Hom_A(C,I_1) \to \Hom_A(C,I_0) \to M \to 0,
    \end{displaymath}
    and also, this sequence stays exact when we apply to it the
    functor $\Hom_A(\Hom_A(C,J),-)$ for any injective $A$--module $J$.
  \end{enumerate}
  \smallskip
  \noindent
  $M$ is called \emph{$C$--Gorenstein projective} if:
  \begin{enumerate}
  \item[(P1)] $\Ext_A^{\geqslant 1}(M,C\ten_AP)=0$ for all projective
    $A$--modules $P$.
  \item[(P2)] There exist projective $A$--modules $P^0,P^1,\ldots$
    together with an exact sequence:
    \begin{displaymath}
       0 \to M \to C\ten_AP^0 \to C\ten_AP^1 \to \cdots,
    \end{displaymath}
    and furthermore, this sequence stays exact when we apply to it the
    functor \mbox{$\Hom_A(-,C\ten_AQ)$} for any projective $A$--module
    $Q$. 
  \end{enumerate}
  \smallskip
  \noindent
  Finally, $M$ is called \emph{$C$--Gorenstein flat} if:
  \begin{enumerate}
  \item[(F1)] $\Tor^A_{\geqslant 1}(\Hom_A(C,I),M)=0$ for all
    injective $A$--modules $I$.
  \item[(F2)] There exist flat $A$--modules $F^0,F^1,\ldots$ together
    with an exact sequence:
    \begin{displaymath}
       0 \to M \to C\ten_AF^0 \to C\ten_AF^1 \to \cdots,
    \end{displaymath}
    and furthermore, this sequence stays exact when we apply to it the
    functor $\Hom_A(C,I)\ten_A-$ for any injective $A$--module $I$.
  \end{enumerate}
\end{Definition}

\begin{Example} \label{exa:injective_is_CGinjective}
  (a) If $I$ is an injective $A$--module, then $\Hom_A(C,I)$ and $I$
  are $C$--Gorenstein injective because:

  \smallskip
  \noindent
  It is easy to see that $\Hom_A(C,I)$ is $C$--Gorenstein injective.
  Concerning $I$ itself it is clear that condition (I1) of Definition
  \ref{def:Cmodules} is satisfied. From Lemma
  \ref{lem:Proj/A=>GorProj/AxC} it follows that $I$ is Gorenstein
  injective over $A \skewtimes C$, so iterating Lemma
  \ref{lem:CGorInj/A<=GorInj/AxC} we also get condition (I2).

  \medskip
  \noindent
  (b) Similarly, if $P$ is a projective $A$--module, then $C\ten_AP$
  and $P$ are $C$--Gorenstein projective. The last claim uses Lemmas
  \ref{lem:Proj/A=>GorProj/AxC} and \ref{lem:CGorProj/A<=GorProj/AxC}.

  \medskip
  \noindent
  (c) If $F$ is a flat $A$--module, then $C\ten_AF$ and $F$ are
  $C$--Gorenstein flat. The last claim uses (a) together with
  Propositions \ref{prop:dual}, \ref{prp:inj_proj_module_case}(1),
  \ref{prp:flat_module_case} (the last two can be found below).
\end{Example}

\begin{Definition} \label{def:new_dimensions}
  By Example \ref{exa:injective_is_CGinjective}(a), there exists for
  every homologically left-bounded complex $N$ a left-bounded complex
  $Y$ of $C$--Goren\-stein injective modules with $Y \simeq N$ in
  $\sD(A)$ (as one could take $Y$ to be an injective resolution of
  $N$). Every such $Y$ is called a \emph{$C$--Gorenstein injective
    resolution} of $N$.

  \medskip
  \noindent
  \emph{$C$--Gorenstein projective} and \emph{$C$--Gorenstein flat
    resolutions} of homologically right-bounded complexes are defined
  in a similar way, and they always exist by Examples
  \ref{exa:injective_is_CGinjective}(b) and (c). Thus, we may define:

  \medskip
  \noindent
  For any homologically left-bounded $A$--complex $N$ we introduce:
  \begin{displaymath}
    \GIdim N \,=\, \underset{Y}{\inf} \Big( \sup
    \big\{n\in\BZ \,|\, Y_{-n} \neq 0 \big\} \Big), 
  \end{displaymath}
  where the infimum is taken over all $C$--Gorenstein injective
  resolutions $Y$ of $N$. For a homologically right-bounded
  $A$--complex $M$ we define:
  \begin{displaymath}
    \GPdim M \,=\, \underset{X}{\inf} \Big( \sup
    \big\{n\in\BZ \,|\, X_n \neq 0 \big\} \Big), 
  \end{displaymath}
  where the infimum is taken over all $C$--Gorenstein projective
  resolutions $X$ of $M$. Finally, we define $\GFdim M$ anologously
  to $\GPdim M$. 
\end{Definition}

\begin{Observation}
  Note that when $C=A$ in Definition \ref{def:Cmodules}, we recover
  the categories of ordinary Gorenstein injective, Gorenstein
  projective, and Gorenstein flat $A$--modules.
 
  \medskip
  \noindent  
  Thus, $\AGIdim(-)$, $\AGPdim(-)$, and $\AGFdim(-)$ are the usual
  Gorenstein injective, Gorenstein projective and Gorenstein flat
  dimensions over $A$, which one usually denotes $\Gid_A(-)$,
  $\Gpd_A(-)$ and $\Gfd_A(-)$, respectively.
\end{Observation}

\begin{Lemma} \label{lem:CGorInj/A=>GorInj/AxC}
  Let $M$ be an $A$--module which is $C$--Gorenstein injective. Then
  there exists a short exact sequence of $(A \skewtimes C)$--modules,
  \begin{displaymath}
     0 \rightarrow M^{\prime} \longrightarrow U \longrightarrow M 
    \rightarrow 0, 
  \end{displaymath}
  where $U$ is injective over $A \skewtimes C$ and $M^{\prime}$ is
  $C$--Gorenstein injective over $A$.  Furthermore, the sequence stays
  exact if one applies to it the functor $\Hom_{A \skewtimes C}(V,-)$
  for any injective $(A \skewtimes C)$--module $V$.
\end{Lemma}

\begin{proof}
  Since $M$ is $C$--Gorenstein injective, we in particular get a short
  exact sequence of $A$--modules:
  \begin{displaymath}
     0 \rightarrow N \longrightarrow \Hom_A(C,I) \longrightarrow M
     \rightarrow 0,
  \end{displaymath}
  where $I$ is injective and $N$ is $C$--Gorenstein injective, which
  stays exact under $\Hom_A(\Hom_A(C,J),-)$ when $J$ is injective.
  Applying the functor $\Hom_A(-,I)$ to the exact sequence:
  \begin{displaymath} 
    \tag{\text{$*$}}
    0 \rightarrow C \longrightarrow A \skewtimes C \longrightarrow A
    \rightarrow 0
  \end{displaymath}
  gives an exact sequence of $(A \skewtimes C)$--modules:
  \begin{displaymath} 
    \tag{\text{$**$}}
    0 \rightarrow I \longrightarrow \Hom_A(A \skewtimes C,I)
    \longrightarrow \Hom_A(C,I) \rightarrow 0. 
  \end{displaymath}
  If viewed as a sequence of $A$--modules then this is split, because
  the same holds for $(*)$. Combining these data gives a commutative
  diagram of $(A \skewtimes C)$--modules with exact rows:
  \begin{displaymath}
    \xymatrix{0 \ar[r] & M' \ar[r] \ar[d] & \Hom_A(A \ltimes
      C,I) \ar[r] \ar@{->>}[d] & M \ar[r] \ar@{=}[d] & 0 \\
      0 \ar[r] & N \ar[r] & \Hom_A(C,I) \ar[r] & M \ar[r] & 0. } 
  \end{displaymath}
  We will prove that the upper row here has the properties claimed in
  the lemma:
  
  \smallskip
  \noindent
  First, $\Hom_A(A \skewtimes C,I)$ is an injective $(A \skewtimes
  C)$--module by Lemma \ref{lem:induced_injectives}(1). Secondly,
  using the Snake Lemma on the diagram embeds the vertical arrows into
  exact sequences. The leftmost of these is:
  \begin{displaymath}
    0 \rightarrow I \longrightarrow M^{\prime} \longrightarrow N
    \rightarrow 0, 
  \end{displaymath}
  proving that as $A$--modules, $M^{\prime} \cong I \oplus N$.  Here
  $N$ is $C$--Gorenstein injective by construction, and $I$ is by
  Example \ref{exa:injective_is_CGinjective}(a). So $M^{\prime}$ is
  clearly also $C$--Gorenstein injective.

  \medskip
  \noindent
  Finally, by construction, the lower row in the diagram stays exact
  under $\Hom_A(\Hom_A(C,J),-)$ when $J$ is injective. If viewed as a
  sequence of $A$--modules then the sequence $(**)$ is split, so the
  surjection $\Hom_A(A \skewtimes C,I) \longrightarrow \Hom_A(C,I)$ is
  split, and therefore the upper row in the diagram also stays exact
  under $\Hom_A(\Hom_A(C,J),-)$.

  \medskip
  \noindent
  By aplying $\H_0(-)$ to Lemma \ref{lem:calculation}(1), we see that
  the upper row in the diagram stays exact under $\Hom_{A \skewtimes
    C}(\Hom_A(A \skewtimes C,J),-)$ when $J$ is an injective
  $A$--module.  Thus, it also stays exact under $\Hom_{A \skewtimes
    C}(V,-)$ for any injective $(A \skewtimes C)$--module $V$, because
  of Lemma \ref{lem:induced_injectives}(2).
\end{proof}

\noindent
By a similar argument we get:

\begin{Lemma} \label{lem:CGorProj/A=>GorProj/AxC}
  Let $M$ be an $A$--module which is $C$--Gorenstein projective.  Then
  there exists a short exact sequence of $(A \skewtimes C)$--modules,
  \begin{displaymath}
     0 \rightarrow M \longrightarrow R \longrightarrow M^{\prime} 
    \rightarrow 0, 
  \end{displaymath}
  where $R$ is projective over $A \skewtimes C$ and $M^{\prime}$ is
  $C$--Gorenstein projective over $A$.  Furthermore, the sequence stays
  exact if one applies to it the functor $\Hom_{A \skewtimes C}(-,S)$
  for any projective $(A \skewtimes C)$--module $S$.
\end{Lemma}

\begin{Proposition} \label{prp:inj_proj_module_case}
  For any $A$--module $M$ the two conclusions hold:
  \begin{enumerate}
  \item $M$ is $C$--Gorenstein injective if and only if $M$ is
    Gorenstein injective over $A \ltimes C$.
  \item $M$ is $C$--Gorenstein projective if and only if $M$ is
    Gorenstein projective over $A \ltimes C$.
  \end{enumerate}
\end{Proposition}

\begin{proof}
  (1) If $M$ is $C$--Gorenstein injective, then Lemma
  \ref{lem:CGorInj/A=>GorInj/AxC} gives the ``left half'' of a
  complete injective resolution of $M$ over $A \ltimes C$.

  \smallskip
  \noindent
  Conversely, if $M$ is Gorenstein injective over $A \ltimes C$, then
  Lemma \ref{lem:CGorInj/A<=GorInj/AxC} gives the existence of a
  sequence like the one in Definition \ref{def:Cmodules} (I2).  Now,
  to finish the proof we only need to refer to Corollary
  \ref{cor:Ext}(1).

  \smallskip
  \noindent
  (2) Similar, but using Lemmas \ref{lem:CGorProj/A=>GorProj/AxC},
  \ref{lem:CGorProj/A<=GorProj/AxC} and Corollary \ref{cor:Ext}(2).
\end{proof}

\noindent
Before turning to $C$--Gorenstein flat modules, we need to recall
the notion of \emph{Kaplansky classes} from \cite[Definition
2.1]{EnochsLopezRamos}, which is reformulated in Definition
\ref{def:Kaplansky}, Section \ref{sec:proper_dimensions}. The
following lemma will be central:

\begin{Lemma} \label{lem:C_tensor_flad_Kaplansky}
  The class ${\sf F} = \{ C\ten_AF \,|\ F \,\text{flat }
  A\text{--module}\}$ is Kaplansky, and furthermore it is closed under
  direct limits.
\end{Lemma}

\begin{proof}
  Every homomorphism $\varphi \colon C\ten_AF_1 \to C\ten_AF_2$, where
  $F_i$ is flat, has the form $\varphi = C\ten_A\psi$ for some
  homomorphism $\psi \colon F_1 \to F_2$; namely $\psi =
  \Hom_A(C,\varphi)$, because $\Hom_A(C,C\ten_AF_i) \cong F_i$. 

  \medskip
  \noindent
  With this observation in mind it is clear that $\sF$ is closed under
  direct limits, since the class of flat modules has this property.

  \medskip
  \noindent
  To see that $\sF$ is Kaplansky, we first note that a finitely
  generated $A$--module has cardinality at most $\kappa =
  \max\{|A|,\aleph_0\}$.

  \medskip
  \noindent
  Now, assume that $x$ is an element of $G = C\ten_AF$, where $F$ is a
  flat $A$--module. Write $x=\sum_{i=1}^n c_i\ten x_i$ for some
  $c_1,\ldots,c_n \in C$ and $x_1,\ldots,x_n \in F$.  Let $S$ be the
  $A$--submodule of $F$ generated by $x_1,\ldots,x_n$, and then use
  \cite[Lemma 2.5.2]{Xu} (or \cite[Lemma 5.3.12]{EnochsJendabook}) to
  enlarge $S$ to a pure submodule $F'$ in $F$ with cardinality:
  \begin{displaymath}
    |F'| \,\leqslant\, \max\{|S|\!\cdot\!|A|,\aleph_0\} \,\leqslant\,
    \kappa.  
  \end{displaymath}
  Since $F$ is flat and $F' \subseteq F$ is a pure submodule, then
  $F'$ and $F/F'$ are flat as well. Furthermore, exactness of:
  \begin{displaymath} 
    0 \to C\ten_AF' \to C\ten_AF \to C\ten_A(F/F') \to 0
  \end{displaymath}   
  shows that $G' = C\ten_AF'$ is a submodule of $G = C\ten_AF$ which
  contains $x$. Clearly, $G'$ and $G/G' \cong C\ten_A(F/F')$ belong
  to ${\sf F}$, and:
  \begin{displaymath} 
    |G'| \,=\, |C\ten_AF'| \,\leqslant\, |\mathbb{Z}^{(C \times F')}|
    \,\leqslant\, |(2^{\mathbb{Z}})^{(C
      \times F')}| \,=\, |2^{(\mathbb{Z} \times C \times F')}|
    \,\leqslant\,  2^{\kappa}. 
  \end{displaymath}
  The last inequality comes from the fact that all three cardinal
  numbers $|\mathbb{Z}|$, $|C|$ and $|F'|$ are less than $\kappa$.
  Note that the cardinal number $2^{\kappa}$ only depends on the ring
  $A$.
\end{proof}

\noindent
The next proof is modelled on that of \cite[Theorem
(6.4.2)]{Wintherbook}.

\begin{Proposition} \label{prp:flat_module_case}
  Let $M$ be an $A$--module. Then $M$ is $\,C$--Gorenstein flat if and
  only if $M$ is Gorenstein flat over $A \ltimes C$.  In the
  affirmative case, $M$ has the next property, which implies
  Definition \ref{def:Cmodules} \textnormal{(F2)}:
  \begin{enumerate}
  \item[(F2')] There exist flat $A$--modules $F^0,F^1,\ldots$ together
    with an exact sequence:
    \begin{displaymath}
       0 \to M \to C\ten_AF^0 \to C\ten_AF^1 \to \cdots, 
    \end{displaymath}
    and furthermore, this sequence stays exact when we apply to it the
    functor $\Hom_A(-,C\ten_AG)$ for any flat $A$--module $G$. \smallskip
  \end{enumerate}
\end{Proposition}

\begin{proof}
  For the first statement, it suffices by Propositions \ref{prop:dual}
  and \ref{prp:inj_proj_module_case}(1) to show that if $E$ is a
  faithfully injective $A$--module, then:
  \begin{align*}
    \text{$M$ is $C$--Gorenstein flat $\,\Leftrightarrow\,$
      $\Hom_A(M,E)$ is $C$--Gorenstein injective.}
  \end{align*}
  For any injective $A$--module $I$ we have (adjointness)
  isomorphisms:
  \begin{align*}
    \Ext_A^i(\Hom_A(C,I),\Hom_A(M,E)) \,\cong\, \\
    \Hom_A(\Tor^A_i(\Hom_A(C,I),M),E).
  \end{align*}
  Thus, Definition \ref{def:Cmodules} (F1) for $M$ is equivalent to
  (I1) for $\Hom_A(M,E)$. The rest of the proof will concern the
  conditions (F2) for $M$ and (I2) for $\Hom_A(M,E)$ in Definition
  \ref{def:Cmodules}.

  \medskip
  \noindent 
  If $\,\mathbb{S} = 0 \to M \to C\ten_AF^0 \to C\ten_AF^1 \to \cdots$
  is a sequence for $M$ like the one in Definition \ref{def:Cmodules}
  (F2), then, using adjointness, it is easy to see that
  $\Hom_A(\mathbb{S},E)$ is a sequence for $\Hom_A(M,E)$ like the one
  in (I2). Therefore, we have proved the implication ``$\Rightarrow$''

  \medskip
  \noindent
  To show ``$\Leftarrow$'', we assume that $\Hom_A(M,E)$ is
  $C$--Gorenstein injective. As already noted, we only have to show
  Definition \ref{def:Cmodules} (F2) for $M$.  First note that (F2')
  really implies Definition \ref{def:Cmodules} (F2), since:
  \begin{align*}
    \Hom_A(\Hom_A(C,I)\ten_A-,E) &\,\simeq\,
    \Hom_A(-,\Hom_A(\Hom_A(C,I),E)) \\
    &\,\simeq\, \Hom_A(-,C\ten_A\Hom_A(I,E)),
  \end{align*}
  and when $I$ is injective, then $G=\Hom_A(I,E)$ is flat.  In order
  prove (F2'), it suffices to show the existence of a short exact
  sequence:
  \begin{displaymath}
    \tag{\text{$\dagger$}}
    0 \to M \to C\ten_AF \to M' \to 0,
  \end{displaymath}
  satisfying the following three conditions:
  \begin{enumerate}
  \item $F$ is flat, 
  \item $\Hom_A(M',E)$ is $C$--Gorenstein injective, 
  \item \mbox{$\Hom_A((\dagger),C\ten_AG)$} is exact for any flat
    $A$--module $G$.
  \end{enumerate}
  Because then one obtains the sequence in (F2') by iterating
  $(\dagger)$.  By Lemma \ref{lem:C_tensor_flad_Kaplansky}, the class
  of $A$--modules:
  \begin{displaymath}
    {\sf F} \,=\, \big\{ C\ten_AF \,|\ F \text{ flat }
    A\text{--module} \big\}.
  \end{displaymath}
  is Kaplansky. Furthermore, it is closed under arbitrary direct
  products; since $C$ is finitely generated and $A$ is noetherian, and
  hence \cite[Theorem 2.5]{EnochsLopezRamos} implies that every
  $A$--module has an ${\sf F}$--preenvelope.
  
  \medskip
  \noindent
  Note that since $\Hom_A(M,E)$ is $C$--Gorenstein injective, there in
  particular exists an epimorphism $\Hom_A(C,I) \twoheadrightarrow
  \Hom_A(M,E)$, where $I$ is injective.  Applying $\Hom_A(-,E)$, we
  get a monomorphism:
  \begin{align*}
    M &\,\hookrightarrow\, \Hom_A(\Hom_A(M,E),E) \\ 
      &\,\hookrightarrow\, \Hom_A(\Hom_A(C,I),E) \,\cong\,
    C\ten_A\Hom_A(I,E) \,\in\, {\sf F}. 
  \end{align*}
  Thus, $M$ can be embedded into a module from ${\sf F}$. Therefore,
  taking an ${\sf F}$--preenvelope $\varphi \colon M \to C\ten_AF$ of
  $M$, it is automaticly injective; and defining
  $M'=\text{Coker}\,\varphi$, we certainly get an exact sequence
  $(\dagger)$ satisfying (1) and (3).  

  \medskip
  \noindent
  Finally, we argue that (2) is true. Keeping Proposition
  \ref{prp:inj_proj_module_case}(1) in mind we must prove that
  $\Hom_A(M',E)$ is Gorenstein injective over $A \ltimes C$. Applying
  $\Hom_A(-,E)$ to $(\dagger)$ we get:
  \begin{displaymath}
    \tag{\text{$\ddagger$}}
    0 \to \Hom_A(M',E) \to \Hom_A(C,J) \to \Hom_A(M,E) \to 0,
  \end{displaymath}
  where $J \cong \Hom_A(F,E)$ is injective. $\Hom_A(C,J)$ and
  $\Hom_A(M,E)$ are both Gorenstein injective over $A \ltimes C$ ---
  the last module by assumption. Hence, if we can prove that
  $\Ext^1_{A\ltimes C}(U,\Hom_A(M',E))=0$ for every injective $(A
  \ltimes C)$--module $U$, then \cite[Theorem 2.13]{EnochsJenda} gives
  the desired conclusion. Using Corollary \ref{cor:Ext}(1), we must
  prove that:
  \begin{displaymath} 
    \tag{\text{$\natural$}}
    \Ext_A^1(\Hom_A(C,I),\Hom_A(M',E))=0
  \end{displaymath}
  for all injective $A$--modules $I$. Consider the commutative
  diagram with exact columns:
  {\small
  \begin{displaymath}
    \xymatrix{0 & {} \\ 
      \Ext^1_A(\Hom_A(C,I),\Hom_A(M',E)) \ar[u] & 0 \\
      \Hom_A(\Hom_A(C,I),\Hom_A(M,E)) \ar[u] &
      \Hom_A(\Hom_A(C,I)\ten_AM,E) \ar[l]_-{\cong} \ar[u] \\
      \Hom_A(\Hom_A(C,I),\Hom_A(C,J)) \ar[u] &
      \Hom_A(\Hom_A(C,I)\ten_A(C\ten_AF),E) \ar[l]_{\cong} \ar[u] 
    }
  \end{displaymath} }
{\flushleft The} first column is the induced long exact sequence which
comes from applying $\Hom_A(\Hom_A(C,I),-)$ to $(\ddagger)$. We get
another monomorphism when we apply $\Hom_A(C,I)\ten_A-$ to the one $0
\to M \to C\ten_AF$ from $(\dagger)$; this follows from the property
(3) which $(\dagger)$ satisfies together with the calculation
preceding $(\dagger)$.  Turning this into an epimorphism with
$\Hom_A(-,E)$ we get the second column. The vertical isomorphisms are
by adjointness.  The diagram implies that the module in $(\natural)$
is zero.
\end{proof}

\begin{Theorem} \label{thm:main_theorem}
  For any (appropriately homologically bounded) $A$--complex $M$, we
  have the following equalities:
  \begin{eqnarray*}
    \GIdim  M & = & \CGid M, \\ 
    \GPdim  M & = & \CGpd M, \\ 
    \GFdim  M & = & \CGfd M.
  \end{eqnarray*}
\end{Theorem}

\begin{proof}
  The proof uses Propositions \ref{prp:inj_proj_module_case}(1),(2)
  and \ref{prp:flat_module_case} in combination with \cite[Theorems
  (2.5), (2.2) and (2.8)]{CFH}. We only prove that $\GIdim  M = \CGid
  M$, since the proofs of the other two equalities are similar:

  \smallskip
  \noindent
  From Proposition \ref{prp:inj_proj_module_case}(1) we get that every
  $C$--Gorenstein injective $A$--module is also Gorenstein injective
  over $A \ltimes C$, and this give us the inequality ``$\geqslant$''.
 
  \medskip
  \noindent
  For the opposite inequality ``$\leqslant$'', we may assume that
  $n=\CGid M$ is finite. Pick a left-bounded complex $I$ of
  injective $A$--modules such that $I \simeq M$ in $\sD(A)$. By Lemma
  \ref{lem:Proj/A=>GorProj/AxC} the modules $I_i$ are Gorenstein
  injective over $A \ltimes C$, and therefore \cite[Theorem
  (2.5)]{CFH} implies that the $A$--module $\textnormal{Z}^I_{-n}$ is 
  Gorenstein injective over $A \ltimes C$. 

  \medskip
  \noindent
  Now, Proposition \ref{prp:inj_proj_module_case}(1) shows that
  $\textnormal{Z}^I_{-n}$ is $C$--Gorenstein injective. By Example
  \ref{exa:injective_is_CGinjective}(a), the complex
  \mbox{$I_{-n}\!\!\supset \ = \ \cdots \to I_{-n+1} \to
    \textnormal{Z}^I_{-n} \to 0$} consists of $C$--Gorenstein
  injective $A$--modules, and since \mbox{$I_{-n}\!\!\supset \ \simeq
    I \simeq M$} we see that $\GIdim M \leqslant n$.
\end{proof}

\begin{Corollary} \label{cor:C=A}
  For any (appropriately homologically bounded) $A$--complex $M$, we
  have the following equalities:
  \begin{displaymath}
  \begin{array}{lclcl}
    \AGid M & = & \Gid_{A[x]/(x^2)}M & = & \Gid_AM, \\ 
    \AGpd M & = & \Gpd_{A[x]/(x^2)}M & = & \Gpd_AM, \\ 
    \AGfd M & = & \Gfd_{A[x]/(x^2)}M & = & \Gfd_AM.
  \end{array}
  \end{displaymath}
\end{Corollary}

\begin{proof}
  This follows immediately from Theorem \ref{thm:main_theorem}; we
  only have to note that $A \ltimes A \cong A[x]/(x^2)$ (sometimes
  refered to as the \emph{dual numbers} over $A$).
\end{proof}

\noindent
Having realized that, on the level of $A$--complexes, the three
(classical) Gorenstein dimensions can not distinguish between $A$ and
$A \ltimes A$, we can reap a nice result from the work of
\cite{HHPJ2}:

\begin{Theorem} \label{thm:GorThm}
  If $(A,\mathfrak{m},k)$ is local, then the following conditions are
  equivalent: 
  \begin{enumerate}
  \item $A$ is Gorenstein.
  \item There exists an $A$--complex $M$ such that all three numbers
    $\,\fd_AM$, $\Gid_AM$ and $\textnormal{width}_AM$ are finite.
  \item There exists an $A$--complex $N$ such that all three numbers
    $\,\id_AN$, $\Gpd_AN$ and $\textnormal{depth}_AN$ are finite.
  \item There exists an $A$--complex $N$ such that all three numbers
    $\,\id_AN$, $\Gfd_AN$ and $\textnormal{depth}_AN$ are finite.
  \end{enumerate}
\end{Theorem}

\begin{proof}
  It is well-known that over a Gorenstein ring, every homologically
  bounded complex has finite Gorenstein injective, Gorenstein
  projective and Gorenstein flat dimension, and thus $(1) \Rightarrow
  (2), (3), (4)$.

  \medskip
  \noindent
  Of course, $(3) \Rightarrow (4)$; and using Corollary \ref{cor:C=A},
  the remaining implications $(2) \Rightarrow (1)$ and $(4)
  \Rightarrow (1)$ follow immediately from \cite[Propositions 4.5 and
  4.7]{HHPJ2}.
\end{proof}

\begin{Remark}
  There already exist special cases of this result in the litterature:
  If $A$ admits a dualizing complex, then \cite[(3.3.5)]{Wintherbook}
  compared with \cite[Theorems (4.3) and (4.5)]{CFH} gives Theorem
  \ref{thm:GorThm}. If one drops the assumption that a dualizing
  complex should exists, then Theorem \ref{thm:GorThm} is proved in
  \cite[Corollary (3.3)]{HHRWFGID}, but only for modules.
\end{Remark}


\section{Comparison with Christensen's $\GdimC(-)$}
\label{sec:Christensen}

\noindent
In \cite[Definition (3.11)]{Wintherpaper}, Christensen introduced the
number $\GdimC Z$ for any semi-dualizing complex $C$, and any complex
$Z$ with bounded and finitely generated homology.  When $C = A$ (and
$Z$ is a module), we recover Auslander--Bridger's $G$--dimension by
\cite[Theorem (2.2.3)]{Wintherbook}.

\begin{Proposition} \label{pro:GdimC}
  If $\,C$ is a semi-dualizing $A$--module, and $M$ an $A$--complex
  with bounded and finitely generated homology, then:
  \begin{displaymath}
    \GPdim  M \,=\, \GdimC M.
  \end{displaymath}
\end{Proposition}

\begin{proof}
  By Theorem \ref{thm:main_theorem}, the proposition amounts to:
  \begin{displaymath}
    \tag{\text{$*$}}
    \Gpd_{A \skewtimes C} M \,=\, \GdimC M.
  \end{displaymath}
  The homology of $M$ is bounded and finitely generated over $A$, and
  hence it is also bounded and finitely generated over $A \skewtimes
  C$. So by e.g.~\cite[Theorem (2.12)(b)]{CFH} or \cite[Theorem
  (4.2.6)]{Wintherbook}, the left hand side in $(*)$ equals $\GdimAxC
  M$ (Auslander--Bridger's $G$--dimension over the ring $A \ltimes
  C$). We must therefore prove that:
  \begin{displaymath}
    \tag{\text{$**$}}
    \GdimAxC M \,=\, \GdimC M. 
  \end{displaymath}
  The left hand side is finite precisely if the biduality morphism:
  \begin{displaymath}
    M \longrightarrow \RHom_{A \skewtimes C}(
    \RHom_{A \skewtimes C}(M,A \skewtimes C),A \skewtimes C)
  \end{displaymath}                                       
  is an isomorphism, and the right hand side is finite precisely when
  \begin{displaymath}
    M \longrightarrow \RHom_A(\RHom_A(M,C),C)
  \end{displaymath}
  is an isomorphism. But these two morphisms are equal by Lemma
  \ref{lem:extension_formulae}(3), so the left hand side and right
  hand side of $(**)$ are simultaneously finite.  When the left hand
  side of $(**)$ is finite, it equals:
  \begin{displaymath}
     -\inf \RHom_{A \skewtimes C}(M,A \skewtimes C),
  \end{displaymath}
  and when the right hand side is finite, it is equal to:
  \begin{displaymath}
    -\inf \RHom_A(M,C)
  \end{displaymath}
  But these two numbers are equal by Lemma
  \ref{lem:extension_formulae}(2).
\end{proof}

\begin{Observation}
  Christensen's $\GdimC(-)$ only works when the argument has bounded
  and finitely generated homology, but it has the advantage that $C$
  is allowed to be a semi-dualizing \emph{complex}.

  \medskip
  \noindent
  By Theorem \ref{thm:main_theorem}, we get that for $A$--complexes
  $M$, the $C$--Gorenstein projective dimension $\GPdim M$ agrees with
  the ``ring changed'' Gorenstein projective dimension $\CGpd M$.

  \medskip
  \noindent
  It is not immediately clear how one should make either of these
  dimensions work when $C$ is a semi-dualizing \emph{complex}. Because
  in this case, $A \ltimes C$ becomes a differential graded algebra,
  and the $C$--Gorenstein projective objects in Definition
  \ref{def:Cmodules} (from which we build our resolutions) become
  complexes.

  \medskip
  \noindent
  In \cite[Page 28]{Wintherthesis} we find an interesting comment,
  which makes it even more clear why we run into trouble when $C$ is a
  complex:
  
  \medskip 
  \noindent 
  ``\textsl{On the other hand, let $C$ be a semi-dualizing complex
    with $\textnormal{amp}\,C = s >0$. We are free to assume that
    $\inf C = 0$, and it is then immediate from the definition that
    $\GdimC C = 0$; but a resolution of $C$ must have length at least
    $s$, so the $G$--dimension with respect to $C$ can not be
    interpreted in terms of resolutions.}''

  \medskip
  \noindent
  It is notable that the number $\Gpd_A \RHom_A(C,N)$, $N \in
  \sB_C(A)$, occuring in Theorem \ref{thm:CGpd_and_BC} below makes
  perfect sense even if $C$ is a complex.
\end{Observation}


\section{Interpretations via Auslander and Bass categories}
\label{sec:AB}

\noindent
In this section, we interpret the $C$--Gorenstein homological
dimensions from Section \ref{sec:CGdims} in terms of Auslander and
Bass categories.

\begin{Remark} \label{rmk:AB}
  Let $C$ be a semi-dualizing $A$--complex. In \cite[Section
  4]{Wintherpaper} is considered the adjoint pair of functors:
  \begin{displaymath}
    \xymatrix{ \sD(A) \ar@<0.5ex>[rrr]^-{C\LTensor_A-} & {} & {} &
      \sD(A) \ar@<0.5ex>[lll]^-{\RHom_A(C,-)} }
  \end{displaymath}
  and the full subcategories (where $\sD_{\b}(A)$ is the full
  subcategory of $\sD(A)$ consisting of homologically bounded
  complexes): {\small
  \begin{equation*}
    \sA_C(A) = \left\{ M \in \sD(A) \:
    \left|
    \begin{array}{l}
      \mbox{$M$ and $C \LTensor_A M$ are in $\sD_{\b}(A)$ and} \\
      \mbox{$M \rightarrow \RHom_A(C,C \LTensor_A M)$ is an
        isomorphism} 
    \end{array}
    \right. \!\!
    \right\}
  \end{equation*}
  and
   \begin{equation*}
    \sB_C(A) = \left\{ N \in \sD(A) \: 
    \left|
    \begin{array}{l}
      \mbox{$N$ and $\RHom_A(C,N)$ are in $\sD_{\b}(A)$ and} \\
      \mbox{$C \LTensor_A \RHom_A(C,N) \rightarrow N$ is an
        isomorphism} 
    \end{array}
    \right. \!\!
    \right\}.
  \end{equation*} } 
{\flushleft It} is an exercise in adjoint functors that the adjoint
pair above restricts to a pair of quasi-inverse equivalences of
categories:
  \begin{displaymath}
    \xymatrix{ \sA_C(A) \ar@<0.5ex>[rrr]^-{C\LTensor_A-} & {} & {} &
      \sB_C(A). \ar@<0.5ex>[lll]^-{\RHom_A(C,-)} }
  \end{displaymath}
\end{Remark}

\begin{Theorem} \label{thm:CGid_and_AC}
  For any complex $M \in \sA_C(A)$ we have an equality:
  \begin{displaymath}
    \GIdim M  \,=\, \Gid_A(C\LTensor_AM).
  \end{displaymath}
\end{Theorem}

\begin{proof}
  Throughout the proof we make use of the nice desciptions of the
  \emph{modules} in $\sA_C(A)$ and $\sB_C(A)$ from \cite[Observation
  (4.10)]{Wintherpaper}.

  \medskip
  \noindent
  \textsc{Step} 1: In order to prove the equality $\GIdim M =
  \Gid_A(C\LTensor_AM)$, we first justify the (necessary)
  bi-implication:
  \begin{align*} 
    \tag{\text{$\natural$}} 
    &M \textnormal{ is } C\textnormal{--Gorenstein injective} \quad
    \iff \quad \\ &C\ten_AM \textnormal{ is Gorenstein injective}
  \end{align*}
  for any \emph{module} $M \in \sA_C(A)$.
  
  \medskip
  \noindent
  ``$\Rightarrow$'': By Definition \ref{def:Cmodules}(I2) there is an
  exact sequence:
  \begin{displaymath}
     \tag{\text{$*$}}
    \cdots \to \Hom_A(C,I_1) \to \Hom_A(C,I_0) \to M \to 0,
  \end{displaymath}
  where $I_0,I_1,\ldots$ are injective $A$--modules. Furthermore, we
  have exactness of $\Hom_A(\Hom_A(C,J),(*))$ for all injective
  $A$--modules $J$. 
  
  \medskip
  \noindent
  $M$ belongs to $\sA_C(A)$, and so does $\Hom_A(C,I)$ for any
  injective $A$--module $I$, since $I \in \sB_C(A)$ by
  \cite[Proposition (4.4)]{Wintherpaper}. In particular, $C$ is
  Tor-independent with both of the modules $M$ and $\Hom_A(C,I)$ (two
  $A$--modules $U$ and $V$ are Tor-independent if $\Tor_{\geqslant
    1}^A(U,V)=0$).  Hence the sequence $(*)$ stays exact if we apply
  to it the functor $C\ten_A-$, and doing so we obtain:
  \begin{displaymath}
    \tag{\text{$**$}}
    \cdots \to I_1 \to I_0 \to C\ten_AM \to 0.
  \end{displaymath}
  By similar arguments we see that if we apply $\Hom_A(C,-)$ to the
  sequence $(**)$, then we get $(*)$ back. If $J$ is any injective
  $A$--module, then we have exactness of $\Hom_A(J,(**))$ because:
  \begin{align*}
    \Hom_A(J,(**)) &\,\cong\, \Hom_A(C\ten_A\Hom_A(C,J),(**)) \\
    &\,\cong\, \Hom_A(\Hom_A(C,J),\Hom_A(C,(**))) \\
    &\,\cong\, \Hom_A(\Hom_A(C,J),(*)).
  \end{align*}
  Thus, $(**)$ is a ``left half'' of a complete injective resolution
  of the $A$--module $C\ten_AM$. We also claim that
  $\Ext^i_A(J,C\ten_AM)=0$ for all $i>0$ and all injective
  $A$--modules $J$. First note that:
  \begin{align*}
    \tag{\text{$\diamond$}}
    \Ext^i_A(J,C\ten_AM) &\,\stackrel{\rm (a)}{=}\,
    \H^i\RHom_A(C\LTensor_A\RHom_A(C,J),C\LTensor_AM) \\
    &\,\stackrel{\rm (b)}{\cong}\,
    \H^i\RHom_A(\RHom_A(C,J),\RHom_A(C,C\LTensor_AM)) \\ 
    &\,\stackrel{\rm (c)}{\cong}\, \H^i\RHom_A(\RHom_A(C,J),M) \\
    &\,\cong\, \Ext_A^i(\Hom_A(C,J),M).
  \end{align*}
  Here (a) is follows as $J \in \sB_C(A)$ by \cite[Proposition
  (4.4)]{Wintherpaper}; (b) is by adjointness; and (c) is because $M
  \in \sA_C(A)$.  This last module is zero because $M$ is
  $C$--Gorenstein injective.  These considerations prove that
  $C\ten_AM$ is Gorenstein injective over $A$.

  \medskip
  \noindent
  ``$\Leftarrow$'': If $C\ten_AM$ is Gorenstein injective over $A$,
  we have by definition an exact sequence:
   \begin{displaymath}
     \tag{\text{$\dagger$}}
    \cdots \to I_1 \to I_0 \to C\ten_AM \to 0,
  \end{displaymath}
  where $I_0,I_1,\ldots$ are injective $A$--modules. Furthermore, we
  have exactness of $\Hom_A(J,(\dagger))$ for all injective
  $A$--modules $J$.

  \medskip
  \noindent
  Since $I_0,I_1,\ldots$ and $C\ten_AM$ are modules from $\sB_C(A)$,
  then so are all the kernels in $(\dagger)$, as $\sB_C(A)$ is a
  triangulated subcategory of $\sD(A)$. If $N \in \sB_C(A)$, then
  $\Ext^{\geqslant 1}_A(C,N)=0$, and consequently, the sequence
  $(\dagger)$ stays exact if we apply to it the functor $\Hom_A(C,-)$.
  Doing so we obtain:
  \begin{displaymath}
     \tag{\text{$\ddagger$}}
    \cdots \to \Hom_A(C,I_1) \to \Hom_A(C,I_0) \to M \to 0.
  \end{displaymath}
  If $J$ is any injective $A$--module, then we have exactness of the
  complex $\Hom_A(\Hom_A(C,J),(\ddagger))$ because:
  \begin{align*}
    \Hom_A(\Hom_A(C,J),(\ddagger)) &\,\cong\,
    \Hom_A(\Hom_A(C,J),\Hom_A(C,(\dagger))) \\
    &\,\cong\, \Hom_A(C\ten_A\Hom_A(C,J),(\dagger)) \\
    &\,\cong\, \Hom_A(J,(\dagger)).
  \end{align*}
  Furthermore, $(\diamond)$ above gives that:
  \begin{align*}
    \Ext_A^{\geqslant 1}(\Hom_A(C,J),M) \,\cong\, \Ext^{\geqslant
      1}_A(J,C\ten_AM) \,=\, 0,  
  \end{align*}
  for all injective $A$--modules $J$. The last zero is because
  $C\ten_AM$ is Gorenstein injective over $A$. Hence $M$ is
  $C$--Gorenstein injective.

  \medskip
  \noindent
  \textsc{Step} 2: To prove the inequality $\GIdim M \geqslant
  \Gid_A(C\LTensor_AM)$ for any complex $M \in \sA_C(A)$, we may
  assume that $m=\GIdim M= \CGid M$; cf.~Theorem
  \ref{thm:main_theorem}, is finite.  Since $C\LTensor_AM$ is
  homologically bounded, there exists a left-bounded injective
  resolution $I$ of $C\LTensor_AM$, that is, $I \simeq C\LTensor_AM$
  in $\sD(A)$.

  \medskip
  \noindent
  We wish to prove that the $A$--module $\textnormal{Z}^I_{-m}$ is
  Gorenstein injective. Since $M$ belongs to $\sA_C(A)$, we get
  isomorphisms:
  \begin{displaymath}  
    M \,\simeq\, \RHom_A(C,C\LTensor_AM) \,\simeq\, \RHom_A(C,I)
    \,\simeq\, \Hom_A(C,I).
  \end{displaymath}
  Now, $\Hom_A(C,I)$ is a complex of Gorenstein injective $A \ltimes
  C$--modules, and thus the $A$--module $L :=
  \textnormal{Z}_{-m}^{\Hom_A(C,I)}$ is Gorenstein injective over $A
  \ltimes C$ by \cite[Theorem (2.5)]{CFH}. By Proposition
  \ref{prp:inj_proj_module_case}(1), $L$ is also $C$--Gorenstein
  injective. Note that:
  \begin{displaymath}  
    -m \,=\, -\Gid_{A \ltimes C}M \,\leqslant\, \inf M
    \,\stackrel{\rm (a)}{=}\, \inf(C\LTensor_AM) \,=\, \inf I, 
  \end{displaymath}
  where the equality (a) comes from
  \cite[Lemma(4.11)(b)]{Wintherpaper}. Therefore, $0 \to
  \textnormal{Z}^I_{-m} \to I_{-m} \to I_{-m-1}$ is exact, and
  applying the left exact functor $\Hom_A(C,-)$ to this sequence we
  get an isomorphism of $A$--modules:
  \begin{displaymath}  
    \tag{\text{$\flat$}}
    L \,=\, \textnormal{Z}_{-m}^{\Hom_A(C,I)} \,\cong\,
    \Hom_A(C,\textnormal{Z}_{-m}^I). 
  \end{displaymath}
  We have a degreewise split exact sequence of complexes:
  \begin{displaymath} 
    0 \rightarrow \Sigma^{-m}\textnormal{Z}^I_{-m} \longrightarrow
      I_{-m}\!\!\supset \ \longrightarrow I_{-m+1}\!\!\sqsupset 
      \ \longrightarrow 0,
  \end{displaymath}
  where we have used the notation from \cite[Appendix
  (A.1.14)]{Wintherbook} to denote soft and hard truncations.  Since
  \mbox{$I_{-m+1}\!\!\sqsupset\ $} has finite injective dimension it
  belongs to $\sB_C(A)$ by \cite[Proposition (4.4)]{Wintherpaper}, and
  furthermore,
  \begin{displaymath} 
    I_{-m}\!\!\supset \ \, \simeq\, I \,\simeq\,
    C\LTensor_AM \,\in\, \sB_C(A).
  \end{displaymath}
  Thus, the module $\textnormal{Z}^I_{-m}$ is also in $\sB_C(A)$, as
  $\sB_C(A)$ is a triangulated subcategory of $\sD(A)$.  Consequently,
  the module $L$ from $(\flat)$ belongs to $\sA_C(A)$ and has the
  property that $C\otimes_AL \cong \textnormal{Z}^I_{-m}$.  Therefore,
  the implication ``$\Rightarrow$'' in $(\natural)$ gives that
  $\textnormal{Z}^I_{-m}$ is Gorenstein injective over $A$, as
  desired.

  \medskip
  \noindent
  \textsc{Step} 3: To prove the opposite inequality $\GIdim M
  \leqslant \Gid_A(C\LTensor_AM)$ for any complex $M \in \sA_C(A)$, we
  assume that $n= \Gid_A(C\LTensor_AM)$ is finite.  Pick any
  left-bounded injective resolution $I$ of $C\LTensor_AM$.  Then the
  $A$--module $\textnormal{Z}^I_{-n}$ is Gorenstein injective by
  \cite[Theorem (2.5)]{CFH}.

  \medskip
  \noindent
  As in \textsc{Step} 2 we get $M \simeq \Hom_A(C,I)$, and thus
  it suffices to show that the module:
  \begin{displaymath}  
    N \,:=\, \textnormal{Z}_{-n}^{\Hom_A(C,I)} \,\cong\,
    \Hom_A(C,\textnormal{Z}_{-n}^I). 
  \end{displaymath}
  is $C$--Gorenstein injective, because then \mbox{$M \simeq
    \Hom_A(C,I)_{-n}\!\!\supset$} shows that $\GIdim M \leqslant n$.
  As before we get that $N$ is a module in $\sA_C(A)$ with $C\ten_AN
  \cong \textnormal{Z}_{-n}^I$, which this time is Gorenstein
  injective over $A$.  Therefore, the implication ``$\Leftarrow$'' in
  $(\natural)$ gives that $N$ is $C$--Gorenstein injective.
\end{proof}

\noindent
Using Proposition \ref{prp:inj_proj_module_case}(2), a similar
argument gives:

\begin{Theorem} \label{thm:CGpd_and_BC}
  For any complex $N \in \sB_C(A)$ we have an equality:
  \begin{xxalignat}{3}
      &{\hphantom{\qed}}
      & \GPdim N  \,=\, \Gpd_A \RHom_A(C,N).
      & {}
      && \qed
    \end{xxalignat}
\end{Theorem}

\noindent
From Theorems \ref{thm:CGid_and_AC} and \ref{thm:main_theorem}, and
Proposition \ref{prop:dual} we can easily derive:

\begin{Theorem} \label{thm:CGfd_and_BC}
  For any complex $N \in \sB_C(A)$ we have an equality:
  \begin{displaymath}
    \GFdim N  \,=\, \Gfd_A \RHom_A(C,N).
  \end{displaymath}
\end{Theorem}


\begin{proof}
  Let $E$ be a faithfully injective $A$--module. Since $N \in
  \sB_C(A)$ it is easy to see that $\RHom_A(N,E) \simeq \Hom_A(N,E)$
  is in $\sA_C(A)$. Hence:
  \begin{align*}
    \GFdim N &\,\stackrel{\rm (a)}{=}\, \GIdim  \RHom_A(N,E) \\
    &\,\stackrel{\rm (b)}{=}\, \Gid_A\big(C \LTensor_A
    \RHom_A(N,E)\big) \\  
    &\,\stackrel{\rm (c)}{=}\, \Gid_A \RHom_A(\RHom_A(C,N),E) \\ 
    &\,\stackrel{\rm (d)}{=}\, \Gfd_A \RHom_A(C,N).
  \end{align*}
  Here (a) is by Proposition \ref{prop:dual} and Theorem
  \ref{thm:main_theorem}; (b) is by Theorem \ref{thm:CGid_and_AC}; (c)
  is by the isomorphism \cite[(A.4.24)]{Wintherbook}; and finally, (d)
  is by Proposition \ref{prop:dual} and Corollary \ref{cor:C=A}.
\end{proof}


\noindent
In the rest of this section, we assume that $A$ admits a
\emph{dualizing complex} $D^A$; cf.~\cite[Definition (1.1)]{CFH}. The
canonical homomorphism of rings, $A \to A \skewtimes C$, turns $A
\skewtimes C$ into a finitely generated $A$--module, and thus
\begin{displaymath}
  D^{A \skewtimes C} \,=\, \RHom_A(A \skewtimes C,D^A)
\end{displaymath}
is a dualizing complex for $A \skewtimes C$.

\begin{Lemma} \label{lem:extension_formulae_2}
  There is an isomorphism over $A$,
  \begin{displaymath}
     D^{A \skewtimes C} \LTensor_{A \skewtimes C} A \,\cong\,
     \RHom_A(C,D^A). 
  \end{displaymath}
\end{Lemma}

\begin{proof}
  This is a computation:
  \begin{eqnarray*}
    D^{A \skewtimes C} \LTensor_{A \skewtimes C} A
    & = & \RHom_A(A \skewtimes C,D^A) \LTensor_{A \skewtimes C} A \\ 
    & \stackrel{\rm (a)}{\cong}
    & \RHom_A(\RHom_{A \skewtimes C}(A,A \skewtimes C),D^A) \\
    & \stackrel{\rm (b)}{\cong}
    & \RHom_A(C,D^A),
  \end{eqnarray*}
  where (a) holds because $D^A$ has finite injective dimension over
  $A$ and (b) is by Lemma \ref{lem:extension_formulae}(4).
\end{proof}

\noindent
By \cite[Corollary (2.12)]{Wintherpaper}, the complex $C^{\dagger} =
\RHom_A(C,D^A)$ is semi-dualizing for $A$. We now have the following
generalization of the main results in \cite[Theorems (4.3) and
(4.5)]{CFH}:

\begin{Theorem} \label{thm:AB}
  Let $M$ and $N$ be $A$--complexes such that the homology of $M$ is
  right-bounded and the homology of $N$ is left-bounded. Then:
  \begin{enumerate}
  \item $M \in \sA_{C^{\dagger}}(A)$ $\iff$ $\GPdim  M < \infty$
    $\iff$ $\GFdim  M < \infty$.
  \item $N \in \sB_{C^{\dagger}}(A)$ $\iff$ $\GIdim  N < \infty$.
  \end{enumerate}
\end{Theorem}

\begin{proof}
  Recall that $D^{A \skewtimes C} = \RHom_A(A \skewtimes C,D^A)$ is a
  dualizing complex for $A \skewtimes C$. If $M$ is a complex of
  $A$--modules then
  \begin{eqnarray*}
    C^{\dagger} \LTensor_A M
    & = & \RHom_A(C,D^A) \LTensor_A M \\
    & \stackrel{\rm (a)}{\cong}
    & (D^{A \skewtimes C} \LTensor_{A \skewtimes C} A) \LTensor_A M \\ 
    & \cong & D^{A \skewtimes C} \LTensor_{A \skewtimes C} M
  \end{eqnarray*}
  and
  \begin{eqnarray*}
    \RHom_A(C^{\dagger},M)
    & = & \RHom_A(\RHom_A(C,D^A),M) \\
    & \stackrel{\rm (b)}{\cong}
    & \RHom_A(D^{A \skewtimes C} \LTensor_{A \skewtimes C} A,M) \\
    & \stackrel{\rm (c)}{\cong}
    & \RHom_{A \skewtimes C}(D^{A \skewtimes C},\RHom_A(A,M)) \\
    & \cong & \RHom_{A \skewtimes C}(D^{A \skewtimes C},M),
  \end{eqnarray*}
  where (a) and (b) are by Lemma \ref{lem:extension_formulae_2} and
  (c) is by adjunction. So using the adjoint pair:
  \begin{displaymath}
    \xymatrix{ \sD(A) \ar@<0.5ex>[rrr]^-{C^{\dagger}\LTensor_A-} & {} &
      {} & \sD(A) \ar@<0.5ex>[lll]^-{\RHom_A(C^{\dagger},-)} }
  \end{displaymath}
  on complexes of $A$--modules is the same as viewing these complexes
  as complexes of $(A \skewtimes C)$--modules and using the adjoint
  pair:
  \begin{displaymath}
    \xymatrix{ \sD(A \ltimes C) \ar@<0.5ex>[rrrr]^-{D^{A \skewtimes C}
        \LTensor_{A \ltimes C}-} & {} & {} & {} & \sD(A \ltimes C)  
      \ar@<0.5ex>[llll]^-{\RHom_{A \skewtimes C}(D^{A \skewtimes C},-)}
      } 
  \end{displaymath}
  Hence a complex $M$ of $A$--modules is in $\sA_{C^{\dagger}}(A)$ if
  and only if it is in $\sA_{D^{A \skewtimes C}}(A \ltimes C)$ when
  viewed as a complex of $(A \skewtimes C)$--modules.  If $M$ has
  right-bounded homology, this is equivalent both to $\Gpd_{A
    \skewtimes C}M < \infty$ and $\Gfd_{A \skewtimes C}M < \infty$ by
  \cite[Theorem (4.3)]{CFH}, and by Theorem \ref{thm:main_theorem}
  this is the same as $\GPdim M < \infty$ and $\GFdim M < \infty$.

  \medskip
  \noindent
  So part (1) of the theorem follows, and a similar method using
  \cite[Theorem (4.5)]{CFH} deals with part (2).
\end{proof}


\section{Proper dimensions} \label{sec:proper_dimensions}

\noindent
In this section, we define and study the \emph{proper} variants of the
dimensions from Theorem \ref{thm:main_theorem}. The results to follow
depend highly on the work in \cite{EnochsLopezRamos}.

\medskip
\noindent
In Definition \ref{def:new_dimensions} we introduced the
dimensions $\GIdim(-)$, $\GPdim(-)$ and $\GFdim(-)$ for
$A$--complexes. When $M$ is an $A$--module it is not hard to see that
these dimensions specialize to: {\small
\begin{align*}
   \GIdim M &\,=\, \inf \left\{ n \in \mathbb{N}_0 
      \left|
      \begin{array}{l}
        \mbox{$0 \to M \to I^0 \to \cdots \to I^n \to 0$ is exact} \\
        \mbox{and $I^0,\ldots,I^n$ are $C$--Gorenstein injective} 
      \end{array}
      \right. \!\!
      \right\},
\end{align*} }
{\flushleft and} similarly for $\GPdim M$ and $\GFdim M$.

\begin{Definition}
  Let $\sQ$ be a class of $A$--modules (which contains the
  zero-module), and let $M$ be any $A$--module. A \emph{proper left
    $\sQ$--resolution} of $M$ is a complex (not necessarily exact):
  \begin{displaymath}
    \tag{\text{$\dagger$}}
    \cdots \to Q_1 \to Q_0 \to M \to 0,
  \end{displaymath}
  where $Q_0, Q_1, \ldots \in \sQ$ and such that $(\dagger)$ becomes
  exact when we apply to it the functor $\Hom_A(Q,-)$ for every $Q \in
  \sQ$. A \emph{proper right $\sQ$--resolution} of $M$ is a complex
  (not necessarily exact):
  \begin{displaymath}
    \tag{\text{$\ddagger$}}
    0 \to M \to Q^0 \to Q^1 \to \cdots,
  \end{displaymath}
  where $Q^0, Q^1, \ldots \in \sQ$ and such that $(\ddagger)$ becomes
  exact when we apply to it the functor $\Hom_A(-,Q)$ for every $Q \in
  \sQ$.
\end{Definition}

\begin{Definition} \label{dfn:ProperDimensions}
  Let $\sQ$ be a class of $A$--modules, and let $M$ be any
  $A$--module. If $M$ has a proper left $\sQ$--resolution, then we
  define the \emph{proper left $\sQ$--dimension} of $M$ by:
  \begin{displaymath}
    \cL\textnormal{-dim}_{\footnotesize \sQ}M \, = \, \inf \left\{ n
      \in \mathbb{N}_0  
      \left|
      \begin{array}{l}
        \mbox{$0 \to Q_n \to \cdots \to Q_0 \to M \to 0$ is} \\
        \mbox{a proper left $\sQ$--resolution of $M$} 
      \end{array}
      \right. \!\!
      \right\}.
  \end{displaymath}
  Similarly, if $M$ has a proper right $\sQ$--resolution, then we
  define the \emph{proper right $\sQ$--dimension} of $M$ by:
  \begin{displaymath}
    \cR\textnormal{-dim}_{\footnotesize \sQ}M \, = \, \inf \left\{ n
      \in \mathbb{N}_0 
      \left|
      \begin{array}{l}
        \mbox{$0 \to M \to Q^0 \to \cdots \to Q^n \to 0$ is} \\
        \mbox{a proper right $\sQ$--resolution of $M$} 
      \end{array}
      \right. \!\!
      \right\}.
  \end{displaymath}
\end{Definition}

\begin{Definition} \label{dfn:SpecialProperDimensions}
  We use $\sG\sI_C(A)$, $\sG\sP_C(A)$ and $\sG\sF_C(A)$ to denote the
  classes of $C$--Gorenstein injective, $C$--Gorenstein projective and
  $C$--Gorenstein flat $A$--modules, respectively. 

  \medskip
  \noindent
  A proper right $\sG\sI_C(A)$--resolution is called a \emph{proper
    $C$--Gorenstein injective resolution}, and a proper left
  $\sG\sP_C(A)$/$\sG\sF_C(A)$--resolution is called a \emph{proper
    $C$--Gorenstein projective/flat resolution}.

  \medskip
  \noindent
  Finally, we introduce the (more natural) notation:
  \begin{enumerate}
  \item[$\bullet$] $\ProperGIdim(-)$ for the proper right
    $\sG\sI_C(A)$--dimension,
  \item[$\bullet$] $\ProperGPdim(-)$ for the proper left
    $\sG\sP_C(A)$--dimension, 
  \item[$\bullet$] $\ProperGFdim(-)$ for the proper left
    $\sG\sF_C(A)$--dimension,
  \end{enumerate} 
  whenever these dimensions are defined.
\end{Definition}

\noindent
The next definition is taken directly from \cite[Definition
2.1]{EnochsLopezRamos}:

\begin{Definition} \label{def:Kaplansky}
  Let ${\sf F}$ be a class of $A$--modules. Then ${\sf F}$ is called
  \emph{Kaplansky} if there exists a cardinal number $\kappa$ such
  that for every module $M \in {\sf F}$ and every element $x \in M$
  there is a submodule $N \subseteq M$ satisfying $x \in N$ and $N,
  M/N \in {\sf F}$ with $|N| \leqslant \kappa$.
\end{Definition}

\begin{Lemma} \label{lem:C-Gor-inj_Kaplansky}
  The class of $\,C$--Gorenstein injective $A$--modules is Kaplansky. 
\end{Lemma}

\begin{proof}
  The class of Gorenstein injective $(A \ltimes C)$--modules is
  Kaplansky by \cite[Proposition 2.6]{EnochsLopezRamos}. Let $\kappa$
  be a cardinal number which implements the Kaplansky property for
  this class.
  
  \medskip
  \noindent
  Now assume that $M$ is a $C$--Gorenstein injective $A$--module, and
  that $x \in M$ is an element. By Proposition
  \ref{prp:inj_proj_module_case}(1), $M$ is Gorenstein injective over
  $A \ltimes C$, and thus there exists a Gorenstein injective $(A
  \ltimes C)$--submodule $N \subseteq M$ with $x \in N$ and $|N|
  \leqslant \kappa$, and such that the $(A \ltimes C)$--module $M/N$
  is Gorenstein injective.

  \medskip
  \noindent
  Since $M$ is an $A$--module, when we consider it as a module over $A
  \ltimes C$, it is annihilated by the ideal $C \subseteq A \ltimes
  C$. Consequently, the two $(A \ltimes C)$--modules $N$ and $M/N$ are
  also annihilated by $C$. This means that $N$ and $M/N$ really are
  $A$--modules which are viewed as $(A \ltimes C)$--modules. Hence
  Proposition \ref{prp:inj_proj_module_case}(1) implies that $N$ and
  $M/N$ are $C$--Gorenstein injective $A$--modules; and we are done.
\end{proof}

\begin{Theorem} \label{thm:proper_CGid}
  Every $A$--module $M$ has a proper $C$--Gorenstein injective
  resolution, and we have an equality:
  \begin{displaymath}
    \ProperGIdim M \,=\, \GIdim M. 
  \end{displaymath}
\end{Theorem}

\begin{proof}
  By Lemma \ref{lem:C-Gor-inj_Kaplansky} above, the class of
  $C$--Gorenstein injective $A$--modules is Kaplansky, and it is
  obviously also closed under arbitrary direct products. Therefore,
  \cite[Theorem 2.5 and Remark 3]{EnochsLopezRamos} implies that every
  $A$--module admits a proper $C$--Gorenstein injective resolution.

  \medskip
  \noindent
  Every injective $A$--module is also Gorenstein injective by Example
  \ref{exa:injective_is_CGinjective}(a), and hence a proper
  $C$--Gorenstein injective resolution is exact. Consequently, we
  immediately get the inequality:
  \begin{displaymath}
    \ProperGIdim M \,\geqslant\, \GIdim M. 
  \end{displaymath}
  To show the opposite inequality, we may assume that $n = \GIdim M$
  is finite. Let $0 \to M \to E^0 \to E^1 \to \cdots$ be a proper
  $C$--Gorenstein injective resolution of $M$. Defining $D^n =
  \textnormal{Coker}(E^{n-2} \to E^{n-1})$ we get an exact sequence:
  \begin{displaymath}
    0 \to M \to E^0 \to \cdots \to E^{n-1} \to D^n \to 0,
  \end{displaymath}
  which also stays exact when we apply to it the (left exact) functor
  $\Hom_A(-,E)$ for every $C$--Gorenstein injective $A$--module $E$.
  Since $\GIdim M = \CGid M = n$, we get by \cite[Theorem
  2.22]{HHGorensteinHomDim} and Proposition
  \ref{prp:inj_proj_module_case}(1) that $D^n$ is $C$--Gorenstein
  injective, so $\ProperGIdim M \leqslant n$.
\end{proof}

\noindent
Sometimes, nice proper $C$--Gorenstein injective resolutions exist:

\begin{Proposition} \label{prp:nice_CGI-res}
  If $M$ is module in $\sA_C(A)$ such that $n = \GIdim M$ is finite,
  then there exists a proper $C$--Gorenstein injective resolution of
  the form:
  \begin{displaymath}
    \tag{\text{$*$}}
    0 \to M \to H^0 \to \Hom_A(C,I^1) \to \cdots \to \Hom_A(C,I^n)
    \to 0,
  \end{displaymath}
  where $H^0$ is $C$--Gorenstein injective and $I^1, \ldots, I^n$ are
  injective. 
\end{Proposition}

\begin{proof}
  As in the proof of Theorem \ref{thm:CGid_and_AC}, the assumption $M
  \in \sA_C(A)$ gives the existence of an exact sequence of
  $A$--modules:
  \begin{displaymath}
    0 \to M \to \Hom_A(C,J^0) \to \ldots \to \Hom_A(C,J^{n-1}) \to
    D^n \to 0,
  \end{displaymath}
  where $J^0,\ldots,J^{n-1}$ are injective, and $D^n$ is Gorenstein
  injective over $A \ltimes C$. Applying Lemma
  \ref{lem:CGorInj/A<=GorInj/AxC} to $D^n$ we get a commutative
  diagram of $A$--modules with exact rows: 
  {\small
  \begin{displaymath}
    \xymatrix{0 \ar[r] & M \ar@{..>}[d] \ar[r] & \Hom_A(C,J^0)
      \ar@{..>}[d] \ar[r] & \cdots \ar[r] & \Hom_A(C,J^{n-1}) 
      \ar@{..>}[d] \ar[r] & D^n \ar@{=}[d] \ar[r] & 0 \\  
    0 \ar[r] & D^0 \ar[r] & \Hom_A(C,U^0) \ar[r] & \cdots
    \ar[r] & \Hom_A(C,U^{n-1}) \ar[r] & D^n \ar[r] & 0 } 
  \end{displaymath} } 
$\!\!\!\!$ where $U^0,\ldots,U^{n-1}$ are injective and $D^0$ is
$C$--Gorenstein injective. The mapping cone of this chain map is of
course exact, and furthermore, it has $0 \to D^n
\stackrel{=}{\longrightarrow} D^n \to 0$ as a subcomplex.

  \medskip
  \noindent
  Consequently, we get the exact sequence $(*)$, where $I^i = U^{i-1}
  \oplus J^i$ for $i = 1, \ldots, n-1$ together with $I^n = U^{n-1}$
  are injective; and $H^0 = D^0 \oplus \Hom_A(C,J^0)$ is
  $C$--Gorenstein injective.
 
  \medskip
  \noindent
  We claim that the sequence $(*)$ remains exact when we apply to it
  the functor $\Hom_A(-,N)$ for any $C$--Gorenstein injective
  $A$--module $N$ (and this will finish the proof):
  
  \medskip
  \noindent
  Splitting $(*)$ into short exact sequences, we get sequences of the
  form $0 \to X \to Y \to Z \to 0$, where $Z$ has the property that it
  fits into an exact sequence:
  \begin{displaymath}
    0 \to Z \to \Hom_A(C,E^0) \to \Hom_A(C,E^m) \to 0,
  \end{displaymath}
  where $E^0,\ldots,E^m$ are injective. Therefore, it suffices to
  prove that every such module $Z$ satisfies $\Ext_A^1(Z,N)=0$ for all
  $C$--Gorenstein injective modules $N$. But as $\Ext_A^{\geqslant
    1}(\Hom_A(C,E^i),N)=0$ for $i=0,\ldots,m$, this follows easily.
\end{proof}

\noindent
We do not know if every module has a proper $C$--Gorenstein projective
resolution. However, in the case where $A$ admits a dualizing complex
and where $C=A$, then the answer is positive by \cite[Theorem
3.2]{PJgorproj}.

\medskip
\noindent
``Dualizing'' the proof of Theorem \ref{thm:proper_CGid} (except the
first part about existence of proper resolutions) and Proposition
\ref{prp:nice_CGI-res}, we get:

\begin{Theorem} \label{thm:proper_CGpd}
  Assume that $M$ is an $A$--module which has a proper $C$--Gorenstein
  projective resolution. Then we have an equality:
  \begin{xxalignat}{3}
      &{\hphantom{\qed}}
      & \ProperGPdim M \,=\, \GPdim M.  
      & {}
      && \qed
    \end{xxalignat}
\end{Theorem}

\begin{Proposition} \label{prp:nice_CGP-res}
  If $M$ is module in $\sB_C(A)$ such that $n = \GPdim M$ is finite,
  then there exists a proper $C$--Gorenstein projective resolution of
  the form:
  \begin{displaymath}
    0 \to C\ten_AP_n \to \cdots \to C\ten_AP_1 \to G_0 \to M \to 0
  \end{displaymath}
  where $G_0$ is $C$--Gorenstein projective and $P_1, \ldots, P_n$ are
  projective. Furthermore, if $M$ is finitely generated, then $G_0,
  P_1, \ldots, P_n$ may be taken to be finitely generated as well.
  \hspace{\stretch{1}} $\qed$
\end{Proposition}

\noindent
The $C$--Gorenstein flat case is more subtle. We begin with the next:

\begin{Lemma} \label{lem:C-Gor-flat_Kaplansky}
  The class of $C$--Gorenstein flat $A$--modules is Kaplansky, and
  closed under direct limits.
\end{Lemma}

\begin{proof}
  As in the proof of Lemma \ref{lem:C-Gor-flat_Kaplansky}; this time
  using \cite[Proposition 2.10]{EnochsLopezRamos}, we see that the
  class of $C$--Gorenstein flat $A$--modules is Kaplansky.

  \medskip
  \noindent
  By Proposition \ref{prp:flat_module_case}, a module $M$ is
  $C$--Gorenstein flat if and only if $M$ satisfies conditions (F1) in
  Definition \ref{def:Cmodules} and (F2') in Proposition
  \ref{prp:flat_module_case}. Clearly, the condition (F1) is closed
  under direct limits. 

  \medskip
  \noindent
  Concerning condition (F2'), we recall from Lemma
  \ref{lem:C_tensor_flad_Kaplansky} that the class of $A$--modules
  ${\sf F} = \{ C\ten_AF \,|\ F \,\text{flat } A\text{--module}\}$ is
  closed under direct limits. Condition (F2') states that $M$ admits
  an infinite proper right $\sF$--resolution, or in the language of
  \cite{EnochsJendaOyonarte, EnochsLopezRamos}, that $\mu_{\sf F}(M) =
  \infty$. Hence \cite[Theorem 2.4]{EnochsLopezRamos} implies that
  also (F2') is closed under direct limits.
\end{proof}

\begin{Theorem} \label{thm:proper_CGfd}
  Every $A$--module $M$ has a proper $C$--Gorenstein flat resolution,
  and there is an equality:
  \begin{displaymath}
    \ProperGFdim M \,=\, \GFdim M. 
  \end{displaymath}
\end{Theorem}

\begin{proof}
  The class $\sG\sF_C(A)$ of $C$--Gorenstein flat modules contains the
  projective (in fact, flat) modules by Example
  \ref{exa:injective_is_CGinjective}(c), and furthermore, it is closed
  under extensions by \cite[Theorem 3.7]{HHGorensteinHomDim} and
  Proposition \ref{prp:flat_module_case}.
  
  \medskip
  \noindent
  Thus, by Lemma \ref{lem:C-Gor-flat_Kaplansky} above and
  \cite[Theorem 2.9]{EnochsLopezRamos} we conclude that the pair
  $(\sG\sF_C(A),\sG\sF_C(A)^{\perp})$ is a \emph{perfect cotorsion
    theory} according to \cite[Definition 2.2]{EnochsLopezRamos}. In
  particular, every module admits a $C$--Gorenstein flat (pre)cover,
  and hence proper $C$--Gorenstein flat resolutions always exist.

  \medskip 
  \noindent
  The equality $\ProperGFdim M = \GFdim M$ follows as in Theorem
  \ref{thm:proper_CGid}; this time using \cite[Theorem
  3.14]{HHGorensteinHomDim} instead of \cite[Theorem
  2.22]{HHGorensteinHomDim}.
\end{proof}


\end{document}